\documentclass[psamsfonts]{amsart}
\usepackage{mathpazo}

\usepackage[hidelinks]{hyperref}
\usepackage{calc}
\newsavebox\CBox
\newcommand\hcancel[2][0.5pt]{%
  \ifmmode\sbox\CBox{$#2$}\else\sbox\CBox{#2}\fi%
  \makebox[0pt][l]{\usebox\CBox}%
  \rule[0.5\ht\CBox-#1/2]{\wd\CBox}{#1}}

\usepackage{soul}
\usepackage{amsmath,amsthm,amssymb,amscd,amsfonts,amsbsy}
\usepackage{color}
\usepackage{hyperref,url}
\usepackage{float}
\usepackage{hyperref}
\usepackage{enumerate}
\usepackage{enumitem}   
\usepackage{cancel}
\usepackage{mathrsfs}  
\usepackage{mathrsfs}  
\usepackage{amsthm}
\usepackage{amsmath}
\usepackage{xcolor}
\usepackage{amsfonts}
\usepackage{amssymb}
\usepackage{amsthm}
\usepackage{amsmath}
\usepackage{amsfonts}
\usepackage{amssymb}
\usepackage{mathtools}
\usepackage{mathrsfs}  
\usepackage[normalem]{ulem}
\newcommand*{\defeq }{\mathrel{\vcenter{\baselineskip0.5ex \lineskiplimit0pt
                     \hbox{\scriptsize.}\hbox{\scriptsize.}}}%
                     =}

\usepackage[toc,page]{appendix}
\usepackage{tikz-cd}
\theoremstyle{definition}

\theoremstyle{plain}
\newtheorem{theorem}{Theorem}

\theoremstyle{remark}
\newtheorem {mtheorem} {\bf Theorem}

\def\R{\mathbb R}
\def \d {\mathrm{d}}
\def \x {\mathbf{x}}
\def \y {\mathbf{y}}

\title[A note on Vishik's normal form]
{A note on Vishik's normal form}
\author[M. M. Castro, R. M. Martins, D. D. Novaes]
{Matheus M. Castro, Ricardo M. Martins, and Douglas D. Novaes}

\address{Departamento de Matem\'{a}tica, Universidade
Estadual de Campinas, Rua S\'{e}rgio Buarque de Holanda, 651, Cidade
Universit\'{a}ria Zeferino Vaz, 13083--859, Campinas, SP, Brazil}
\email{manzatto.castro@gmail.com, rmiranda@unicamp.br,ddnovaes@unicamp.br}

\allowdisplaybreaks
\begin{document}



\maketitle

\begin{abstract}
The Vishik's Normal Form provides a local smooth conjugation with a linear vector field for smooth vector fields near contacts with a manifold. In the present study, we focus on the analytic case. Our main result ensures that for  analytic vector field and manifold, the conjugation with the Vishik's normal form is also analytic. As an application, we investigate the analyticity of Poincar\'{e} Half Maps defined locally near contacts between analytic vector field and manifold.
\end{abstract}


\section{Introduction and Main Result}
Let $M$ be a $\mathcal{C}^{\omega}$ (resp. $\mathcal{C}^{\infty}$) manifold and denote by $\mathfrak{X}^\omega (M)$ (resp. $\mathfrak{X}^\infty (M)$) the set of $\mathcal{C}^{\omega}$ (resp. infinitely differentiable $\mathcal{C}^{\infty}$) vector fields defined on $M.$ A usual,  $\mathcal{C}^{\omega}$ and $\mathcal{C}^{\infty}$ mean, respectively, analytic and infinitely differentiable.
Also, consider a codimension-1  $\mathcal{C}^{\omega}$ (resp. $\mathcal{C}^{\infty}$) embedded submanifold $\Sigma$ of $M.$ Given $p\in\Sigma,$ there exist a neighborhood $U_p\subset M$ of $p$ and a $\mathcal{C}^{\omega}$ (resp. $\mathcal{C}^{\infty}$) function $h: U_p\to \mathbb{R}$ such that $\Sigma\cap U_p=h^{-1}(0)$ and $\nabla h(\x)\neq0$ for every $\x\in\Sigma\cap U_p.$ 
In the above setting, we say that $p$ is a {\it contact of order $k$} ({\it $k$-contact}, for short) between $X$ and $\Sigma$ if $0$ is a root of multiplicity $k+1$ of $f(t)\defeq h\circ X_t(p),$ where $t\mapsto X_t(p)$ is the trajectory of $X$ starting at $p.$ Equivalently,
$$Xh(p) = X^2h(p) = \ldots = X^{k}h(p) =0,\text{ and } X^{k+1} h(p)\neq 0.$$
Here, $X^n h(p)$ denotes the $n$-th Lie derivative of $h$ in the direction $X$ at $p\in\Sigma,$ which is defined recursively as $Xh(p) = \nabla h(p)\cdot X(p)$ and  $X^{n}h(p) =  \nabla\left( X^ {n-1}h\right)(p)\cdot X(p),$ for $n>1.$  In addition, we say that a $k$-contact $p$ between $X$ and $\Sigma$ is {\it simple} if
$$\left\{\nabla h(p), \nabla Xh(p) , \ldots, \nabla X^{k}h(p)\right\} $$
is linearly independent. In particular, simple $k$-contacts for $k=1$ and $k=2$ are called fold and cusp contacts, respectively. One can see that the definition above does not depend on the function $h$ in the following sense: if $p$ is a simple $k$-contact point with $\Sigma$ and $f:U_p'\rightarrow\R$ is any other function defined in some neighborhood $U_p'\subset M$ of $p$ such that $\Sigma\cap U_p'=f^{-1}(0)$ and $\nabla f(\x)\neq0$ for every $\x\in\Sigma\cap U_p',$ then $$\left\{\nabla f(p), \nabla Xf(p) , \ldots, \nabla X^{k}f(p)\right\} $$
is also linearly independent.

In 1972, Vishik \cite{V72}  studied vector fields near simple $k$-contacts. Assuming that both, vector field and submanifold, were smooth of class $\mathcal{C}^{\infty},$ he provided local linear normal forms for such vector fields through conjugation by $\mathcal{C}^{\infty}$ maps. Applications of the Vishik's Normal Form can be found mainly in the study of bifurcations in piecewise smooth dynamical systems (see, for instance, 
\cite{AJMT19,BCT12,BMT13,JT03,NTZ18,SM02}).

In the present study, we focus on the analytic case. Our main result ensures that for analytic vector field and manifold the conjugation with the Vishik's Normal Form is also analytic. For the sake of completeness, we also state the $\mathcal{C}^{\infty}$ case.  

\begin{mtheorem}\label{thm:vishik} 
Let $M$ be a $\mathcal{C}^{\omega}$ (resp. $\mathcal{C}^{\infty}$) manifold, $\Sigma$ a codimension-1  $\mathcal{C}^{\omega}$ (resp. $\mathcal{C}^{\infty}$) embedded submanifold of $M,$ and $X \in \mathfrak{X}^{\infty}(M)$ $(\text{resp. } X \in \mathfrak{X}^{ \omega}(M)).$ Suppose that $p$ $\in\Sigma$ is a simple $k$-contact between $X$ and $\Sigma$ with $k\leq m-1,$ $m \defeq  \dim(M).$ Then, there exists a $\mathcal{C}^{\omega}$ (resp. $\mathcal{C}^{\infty}$) diffeomorphism  $\psi:U\rightarrow V,$ where $U\subset M$ and $V\subset \R^m$ are respectively neighborhoods of $p$ and $0,$ such that $\psi(\Sigma\cap U) =\{(x_1,\ldots,x_m)\in V:\,x_1=0\},$
$$\psi_*X(x_1,\ldots,x_m) = (x_2,x_3,\ldots,x_{k+1},1, 0 , \ldots , 0),  \ \text{if}\ k<m-1,$$
and
$$\psi_*X(x_1,\ldots,x_m) = (x_2,x_3,\ldots,x_{k+1},1),  \ \text{if}\ k=m-1,$$
for every $(x_1,\ldots,x_m)\in V.$
\end{mtheorem}

Theorem \ref{thm:vishik} is proved in Section \ref{proof:vishik}. Its proof is based on the following result.

\begin{mtheorem}\label{thm:preparation} Let $\varphi$ be a real-valued $\mathcal{C}^{\omega}$ (resp. $\mathcal{C}^{\infty}$)  function defined in a neighborhood of $(0,0)\in\R\times\R^n.$ Assume that
$$\varphi(0,0)=0,\ \frac{\partial \varphi}{\partial t} (0,0)=0,\ldots, \ \frac{\partial^{k-1} \varphi}{\partial t^{k-1}} (0,0)=0, \ \frac{\partial^{k} \varphi}{\partial t^{k}} (0,0)\neq 0. $$
Then, there exist $\mathcal{C}^{\omega}$ (resp. $\mathcal{C}^{\infty}$) real-valued functions $a_1,\ldots,a_n,b ,$ defined in a neighborhood of $(0,0)\in \mathbb{R}\times \mathbb{R}^n,$ such that
$$t^k + \sum_{i=1}^{k-1} t^i a_i(\varphi(t,x),x) = b(\varphi(t,x),x),\ \text{in a neighborhood}\ (0,0). $$
\end{mtheorem}

Theorem \ref{thm:preparation} is proved in Section \ref{sec:proof1}.

{ As mentioned, the $\mathcal{C}^{\infty}$ version of Theorem \ref{thm:vishik} has been previously obtained by Vishik \cite{V72}. Its proof was based in the $\mathcal{C}^{\infty}$ version of the preparation Theorem \ref{thm:preparation}, which was mentioned by Vishik as a consequence of a Corollary of a Malgrange Preparation Theorem  
\cite{Malg3}.  It is worth mentioning that the Vishik's paper \cite{V72} is originally written in Russian and, as far as we known, does not have a published English version, becoming a difficult to access paper. 

Here, in addition to a detailed proof of the original Vishik's result, we check that it is also valid in the analytic context, that is, for analytic vector field and manifold the conjugation with the Vishik's normal form is also analytic. We emphasize that this fact has not been noticed before and represents an improvement of the original result, which can be used, for instance, to ensure the analyticity of Poincar\'{e} Half Maps, as illustrated in the next section. It is worth mentioning that the proof of Theorem \ref{thm:vishik} in the $\mathcal{C}^{\omega}$ context follows the same steps as the $\mathcal{C}^{\infty}$ case performed in \cite{V72}, differing only in the preparation Theorem \ref{thm:preparation}, which will be verified in the $\mathcal{C}^{\omega}$ context as a consequence of a Corollary of a Weierstrass Preparation Theorem \cite[Chapter II]{Naka}.
}

%
%

\subsection{Application: Poincar\'{e} Half Maps}\label{sec:poincare}

Simple contacts between $X$ and $\Sigma$ admit local return maps of $X$ to $\Sigma,$ usually called {\it Poincar\'{e} Half Maps}. The $\mathcal{C}^{\omega}$ version of the Vishik's Normal Form, Theorem \ref{thm:vishik}, can be used to investigate the analyticity of such maps. 

For instance, let $\Sigma$ be a codimension-1 $\mathcal{C}^{\omega}$ embedded submanifold of a $\mathcal{C}^{\omega}$ manifold $M$ and assume that $p\in\Sigma$ is a fold singularity of $X\in\mathfrak{X}^\omega (M)$ with respect to $\Sigma.$ Theorem \ref{thm:vishik} ensures the existence of a $\mathcal{C}^{\omega}$ chart $(\psi,U)$ such that $\psi_*X(x_1,x_2,\ldots,x_m)=(x_2,1,0,\ldots,0),$ $\psi(p)=0,$ and $\Sigma^*\defeq \psi(\Sigma\cap U)=\{x_1=0\},$ here $m=\dim(M).$ The trajectory of $X^*\defeq\psi_*X$ passing through $(x_1,x_2,\ldots,x_m)$ writes 
$$X^*_t(x_1,x_2,\ldots,x_m)=\left(x_1+t\, x_2+\dfrac{t^2}{2}\ ,\ x_2+t\ ,\ x_3\ ,\ldots,x_m\right).$$ Thus, solving $x_1+t\, x_2+\dfrac{t^2}{2}=0$ for $x_1 = 0$ we get $t=0$ or $t=-2\ x_2.$ Therefore, given $(0,x_2,\ldots,x_m)\in\Sigma^*,$ we can define the Poincar\'{e} half map $P:\Sigma^*\rightarrow\Sigma^*$ by $$P(0,x_2,\ldots,x_m)\defeq X^*_{-2\, x_2}(0,x_2,\ldots,x_m)=(0,-x_2,x_3,\ldots,x_m),$$ which is clearly analytic. Going back through the change of coordinates $\psi,$ we get defined a Poincaré half map $$Q\defeq\psi^{-1}\circ P\circ\psi:\Sigma\cap U\rightarrow\Sigma\cap U,$$  which is a composition of analytic maps and, therefore, analytic. In \cite{CGP01}, using the blow-up method through generalized polar coordinates, the above conclusion has been previously achieved for planar vector fields and $\Sigma$ being a straight line.  As performed in \cite{CGP01}, this can be used to define Lyapunov-like constants for studying the center-focus problem in piecewise analytic vector fields around tangential singularities. 

\section{Preparation Theorem and proof of Theorem \ref{thm:preparation}}\label{sec:preparation}

This section is devoted to the proof of Theorem \ref{thm:preparation}, which is based on a ``Formal Preparation Theorem'' (see Theorem \ref{TeoremaNovo}) regarding finitely generated module over the ring of germs of $\mathcal{C}^{\omega}$ and $\mathcal{C}^{\infty}$ functions. Before the statement of Theorem \ref{TeoremaNovo}, in accordance with \cite{Malgrange}, we introduce some algebraic concepts.

\subsection{Modules}\label{sec:module} Let $\mathscr{R}$ be a comutative ring with unity and $A$ an Abelian group (with respect to the operation $+$). We say that $A$ is a $\mathscr{R}$-module if there exits a map from $\mathscr{R}$ to the set of homomorphism of $A,$ denoted by $ra \in A $ with $r \in \mathscr{R}$ and $a \in A,$ for which the following properties are satisfied.
    \begin{align*}
    (r_1+r_2)a=r_1a+r_2 a,&\  \text{ for } \ r_1,r_2\in \mathscr{R}\text{ and }  \text{ for } a \in A,\\ 
    (r_1 r_2)a= r_1(r_2a),&\  \text{ for } \ r_1,r_2\in \mathscr{R}\text{ and }  \text{ for } a \in A, \\
    r(a_1+a_2) = r a_1 + ra_2,&\  \text{ for } \ r_1,r_2\in \mathscr{R}\text{ and }  \text{ for } a \in A,\\
    1\cdot a = a,& \  \text{ for } a \in A.
      \end{align*}
The $\mathscr{R}$-module $A$ is said to be finitely generated over $\mathscr{R}$ if there exists a finite number of elements $a_1,\ldots,a_n \in A$ such that any element $a\in A$ can be written in the form
$$a = \sum_{i=1}^{n}r_ia_i, \ r_i \in \mathscr{R}. $$


Now, let $\mathscr R$ and $\mathscr R'$ be two commutative rings with unity, $\phi: \mathscr R \to \mathscr R'$ be a ring homomorphism  and $A$ a $\mathscr{R}'$-modulo. Then  $A$  admits a natural structure of $\mathscr{R}$-module considering the map $\phi$ by letting $ra := \phi(r) a,$ for $r\in \mathscr R$ and $a\in A$. Under these hypothesis we say that $A$ can be seen as a $\mathscr{R}$-module via $\phi$.

\subsection{Germs of functions}\label{sec:germ} Let  $M$ be a $\mathcal{C}^{\omega}$  (resp. $\mathcal{C}^{\infty})$ manifold and $p \in M.$ We start by defining the ring $\mathscr C^\omega_p(M)$ (resp. $\mathscr C^\infty_p(M)$) of $\mathcal{C}^{\omega}$  (resp. $\mathcal{C}^{\infty})$ germs of real functions defined on $M$ at $p.$

Let $f:U_p\to \mathbb{R}$ and $g:V_p \to \mathbb{R}$ be $\mathcal{C}^{\omega}$  (resp. $\mathcal{C}^{\infty})$ functions, where $U_p$ e $V_p$ are neighborhoods  of $p$ in $M.$ We say that $f\sim_p g$ if there exists $W_p \subset V_p \cap U_p$ of $p,$ such that $\left.f\right|_{W_p} = \left.g\right|_{W_p}.$ As usual, the set of functions in $\mathcal{C}^{\omega}$  (resp. $\mathcal{C}^{\infty})$ which are $\sim_p$ equivalent to $f$ is denoted by $[f]_p,$ which is called germ of $f$ at $p.$ Thus, we denote by $\mathscr{C}_p^\omega(M)$ (resp. $\mathscr{C}_p^\infty (M)$)  the set of all $\mathcal{C}^{\omega}$  (resp. $\mathcal{C}^{\infty})$ function germs at $p.$ 
Clearly, $\mathscr C ^{\omega, \infty}_{p}(M)\cong \mathscr C ^{\omega, \infty}_{p}(U)$ for any  neighborhood $U\subset M$ of $p.$ Notice that $\mathscr{C}_p^\omega(M)$ (resp. $\mathscr{C}_p^\infty (M)$) has a natural ring structure inherited from $\mathbb{R},$
    $$[f]_p + [g]_p \defeq  [g+f]_p, $$
    $$[f]_p \cdot [g]_p \defeq  [g\cdot f]_p, $$
     where $g+f$ and $g \cdot f$ are defined in $\text{dom}f\cap \text{dom}g.$ Also the neutral element and the unit are given, respectively, by $[0]_p$ and $[1]_p.$     

It is straightforward to see that $\mathscr{M}_p^{\omega} \defeq  \{[f]_p\in \mathscr C^\omega_p;\ f(p) = 0\}$ (resp. $\mathscr{M}_p^{\infty} \defeq  \{[f]_p\in \mathscr C^\infty_p;\ f(p) = 0\}$)
 is the unique maximal ideal of the ring $\mathscr C^\omega_p(M)$ (resp. $\mathscr C^\infty_p(M)$). Accordingly, $\mathscr C^\omega_p(M)$ (resp. $\mathscr C^\infty_p(M)$) is called {\it local ring}.
 
Let $M$ and $N$ be $\mathcal{C}^{\omega}$  (resp. $\mathcal{C}^{\infty})$ manifolds and  $\phi: M \to N$ a $\mathcal{C}^{\omega}$  (resp. $\mathcal{C}^{\infty})$ map. The map $\phi$ induces the following ring homomorphism between $\mathscr C ^{\omega, \infty}_{\phi(p)}(N)$ and $\mathscr C ^{\omega, \infty}_{p}(M),$ 
\begin{align*}
    \phi^* : \mathscr C ^{\omega, \infty}_{\phi(p)}(N) &\to \mathscr C ^{\omega, \infty}_{p}(M)\\
    [f]_{\phi(p)} &\mapsto [f\circ \phi]_p.
\end{align*}
In addition, if $\phi$ is a local $\mathcal C^\omega$ (resp. $\mathcal C^\infty$) diffeomorphism then $\phi^*$ is a ring isomorphism.

From now on, for the sake of simplicity, we shall omit the brackets in $[f]_p$ and denote the germ only by $f.$ 

\subsection{Formal Preparation Theorem} \label{sec:general}

Theorem \ref{thm:preparation} is a consequence of the following preparation theorem.
\begin{theorem}\label{TeoremaNovo}
Let $\phi:U\subset \mathbb R^{n+1}\to \mathbb{R}^{n+1}$ be a $\mathcal{C}^{\omega} (resp. \mathcal C^{\infty})$, defined in a neighborhood $U$ of $0$ such that $\phi(0) = 0,$ 
and  consider $f_1, \ldots, f_n\in \mathscr{C}_0^{\omega,\infty}(\R^{n+1}).$ Then,   $\{f_1, \ldots, f_n\}$ generates $\mathscr{C}_0^{\omega,\infty}(\R^{n+1})$ as a  $\mathscr{C}_0^{\omega,\infty}(\mathbb R^{n+1})$-module via $\phi^*$ if, and only if, the set $$\left\{ f_1(x) + \phi^*(\mathscr{M}_0^{\omega,\infty})\mathscr{C}_0^{\omega,\infty}\left(\mathbb{R}^{n+1}\right),\ldots, f_n(x) + \phi^*(\mathscr{M}_0^{\omega,\infty})\mathscr{C}_0^{\omega,\infty}\left(\mathbb{R}^{n+1}\right)\right\}$$
generates the vector space $\mathscr{C}_0^{\omega,\infty}\left(\mathbb R^{n+1}\right)/\left(\phi^*(\mathscr{M}_0^{\omega,\infty})\mathscr{C}_0^{\omega,\infty}\left(\mathbb{R}^{n+1}\right)\right)$ over the field\\ $\mathscr{C}_0^{\omega,\infty}\left(\mathbb R^{n+1}\right)/\mathscr{M}_0^{\omega,\infty} \cong \R.$
\end{theorem}

For a proof of Theorem \ref{TeoremaNovo} in the $\mathcal{C}^{\omega}$ case see the Corollary of Theorem $1$ of \cite[Chapter II]{Naka} in the particular case  $R=R' =\mathscr{C}_0^{\omega}\left(\R^{n+1}\right)$ and $u=\phi^*$.  It is important to notice that $R'/\mathscr{M}(R)R'$ in \cite{Naka} corresponds  to $\mathscr{C}_0^{\omega,\infty}\left(\mathbb R^{n+1}\right)/\left(\phi^*(\mathscr{M}_0^{\omega,\infty})\mathscr{C}_0^{\omega,\infty}\left(\mathbb{R}^{n+1}\right)\right)$ in our case.  We may also refer to \cite{Houzel}. The  $\mathcal{C}^{\infty}$ case of Theorem \ref{TeoremaNovo} follows from
the equivalence between  $(a)'$ and $(\hat{b})'$ in the Corollary of Theorem 1 of \cite{Malg3}. See also 
\cite{Mal4}.

\subsection{Proof of Theorem \ref{thm:preparation}}\label{sec:proof1} 
 
Denote $\x=(x_1,\ldots,x_n).$ First of all, notice that $\mathscr M_0^{\omega,\infty} =(t,\x) \mathscr{C}_0^{\omega,\infty}(\mathbb{R}^{n+1}),$ where $(t,\x) \mathscr{C}_0^{\omega,\infty}(\mathbb{R}^{n+1})$  denotes the smallest ideal of $\mathscr{C}_0^{\omega,\infty}(\mathbb{R}^{n+1})$ that contains the elements $t,x_1,\ldots,x_n.$ Indeed, if $g(t,\x) \in \mathscr M_0^{\omega,\infty},$ then $g(0,0)=0.$ Thus,
 \begin{align*}
     g(t,\x)&= g(t,\x) - g(0,0)\\
     &= \int_0^{1} \frac{\d }{\d s} g(st,s\x) \ \d s\\
     &= \int_0^{1} t \frac{\partial g}{\partial t} (st,s\x) + x_i \sum_{i=1}^n \frac{\partial g}{\partial x_i} (st,s\x) \ \d s\\
     &= t \left(\int_0^{1}  \frac{\partial g}{\partial t} (st,s\x)\d s\right) + \sum_{i=1}^{n}x_i\left(\int_{0}^{1} \frac{\partial g}{\partial x_i} (st,s\x) \d s \right).
 \end{align*}
Therefore, $g\in (t,\x)\mathscr{C}_0^{\omega,\infty}(\mathbb{R}^{n+1})$ and, consequently, $\mathscr M_0^{\omega,\infty} \subset (t,\x) \mathscr{C}_0^{\omega, \infty}(\mathbb{R}^{n+1}).$ The equality follows from the maximality of  $\mathscr M_0^{\omega,\infty}.$
 
Now, define $f(t,\x)\defeq(\varphi(t,\x),\x).$ We shall prove that $\mathscr{C}_0^{\omega, \infty}(\mathbb{R}^{n+1})$ is a finitely generated as  $\mathscr{C}_0^{\omega, \infty}(\mathbb{R}^{n+1})$-module via $f^*$. Indeed, from the remark above, we obtain
 \begin{equation}\label{relation1} f^*\left(\mathscr M_0^{\omega, \infty}\right) =f^*\left((t,x_1,\ldots,x_n) \mathscr{C}_0^{\omega,\infty}(\mathbb{R}^{n+1})\right) = (\varphi(t,x),\x) \mathscr{C}_0^{\omega,\infty}(\mathbb{R}^{n+1}).
 \end{equation}
We claim that $(\varphi(t,\x),\x) \mathscr{C}_0^{\omega, \infty}(\mathbb{R}^{n+1}) = (t^k,\x)\mathscr{C}_0^{\omega, \infty}(\mathbb{R}^{n+1}).$ Indeed, by {\it Taylor's Theorem}, there exists a $\mathcal C^{\omega,\infty}$ real function $r$ defined in a neighborhood of $0\in \R$, with $r(0)\neq0,$ such that  $\varphi(t,0) = t^k r(t).$
Proceeding as above, we conclude that
$$
\varphi(t,x)=t^k r(t)+ \sum_{i=1}^{n}x_i\left(\int_{0}^{1} \frac{\partial \varphi}{\partial x_i} (t,s\x) \d s \right).
$$
Therefore, $ \varphi(t,\x) \in  (t^k,\x)\mathscr{C}_0^{\omega, \infty}(\mathbb{R}^{n+1}).$ In addition, since $r(0)\neq 0$, we also conclude that  $t^k \in (\varphi(t,\x),\x) \mathscr{C}_0^{\omega,\infty}(\mathbb{R}^{n+1}),$ which proves the claim. So, from \eqref{relation1}, we get that
$$\frac{\mathscr{C}_0^{\omega, \infty}(\mathbb{R}^{n+1})}{f^* \left(\mathscr M_0^{\omega, \infty}\right) \mathscr{C}_0^{\omega, \infty}(\mathbb{R}^{n+1})} = 
 \frac{\mathscr{C}_0^{\omega, \infty}(\mathbb{R}^{n+1})}{(t^k,\x)\mathscr{C}_0^{\omega, \infty}(\mathbb{R}^{n+1})},$$
which has a basis given by 
$\left\{[t^i]_0 + f^* \left(\mathscr M_0^{\omega, \infty}\right)\mathscr{C}_0^{\omega,\infty}(\mathbb{R}^{n+1}):\, i=0,1,\ldots k-1 \right\}.$
Thus, by Theorem \ref{TeoremaNovo}, the set $\left\{[t^i]_0:\, i=0,1,\ldots k-1 \right\}$
generates $\mathscr{C}_0^{\omega,\infty}(\mathbb R^{n+1})$ as a $\mathscr{C}_0^{\omega, \infty}(\mathbb{R}^{n+1})$-module via $f^*$.
Therefore, since $t^k \in \mathscr C_0^{\omega, \infty}(\mathbb R^{n+1}),$ there exists germs $a_1,\ldots,a_n,b \in \mathscr{C}^{\omega, \infty}_0(\mathbb{R}^{n+1})$ such that
$$t^k + \sum_{i=1}^{k-1} t^i a_i(\varphi(t,\x),\x) = b(\varphi(t,\x),\x). $$

\section{Proof of Vishik's normal form}\label{sec:vishik}

This section is devoted to the proof of Theorem \ref{thm:vishik}. First, consider $X\in \mathfrak{X}^r(M)$ and let $\varphi:M \to N$ be a $\mathcal{C}^{\infty}$ diffeomorphism. Defining $Y= \varphi_* X$ and $f:N\to \mathbb{R}$ as $f = h\circ \varphi^{-1},$ one can see by induction that
\begin{equation}\label{prop21}
X^n h(x) = Y^nf(\varphi(x) ),\  \text{ for }n \in \mathbb{N}.
\end{equation}
This will be used later on.

\subsection{The planar case}

Before proving Theorem \ref{thm:vishik}, in order to fix ideas and motivate the techniques that will be used, we shall first study the planar case. 

Consider $M=\mathbb {R}^2,$  $X \in \mathfrak{X}^\omega (M)$ (resp. $X \in \mathfrak{X}^\infty (M)$) a planar vector field, and let $h:\mathbb{R}^2\to \mathbb{R}$ be a $\mathcal{C}^{\omega}$  (resp. $\mathcal{C}^{\infty})$ function for which $0$ is a regular value. Denote $\Sigma = h^ {-1}(0)$ and assume that $(0,0)\in\Sigma$ is a fold contact between $X$ and $\Sigma,$ that is, $$X(0,0)\neq (0,0), \ Xh(0,0) = 0\text{ and } X^2h(0,0) \neq 0.$$
Here, we look for a $\mathcal{C}^{\omega}$  (resp. $\mathcal{C}^{\infty})$ diffeomorphism $\psi=(\psi_1,\psi_2): U\subset \mathbb{R}^2\to V\subset \mathbb{R}^2,$ where $U$ and $V$ are neighborhoods of $0,$ such that
$$\psi_*X(x_1,x_2) = (x_2,1)\text{ and } \psi(\Sigma \cap U)=\{(x_1,x_2)\in V:\,x_1=0\}.$$ 
We anticipate that $\psi$ will be given as the composition of three maps,  $\psi=\gamma\circ\beta\circ\alpha.$ In what follows, we shall construct each one of these maps.

Since $X(0,0)\neq (0,0),$ the {\it Tubular Flow Theorem} provides a chart $(\alpha,U_1)$ of $(0,0)$ satisfying $\alpha(0,0)= (0,0)$ and $$Y(x,y)\defeq\alpha_*X(x,y) = (1,0),\ \text{for every} \ (x,y)\in V_1\defeq \alpha(U_1).$$ Defining $f = h\circ \alpha^{-1},$ we have that $\Sigma$ is transformed into $\Sigma_1\defeq\alpha(\Sigma\cap U_1)=f^{-1}(0)$ and, from \ref{prop21},
\begin{equation}\label{planarR1}
\begin{array}{l}
\displaystyle 0=X h(0,0) = Yf(0,0) =\frac{\partial f}{\partial x}(0,0)\text{ and }\vspace{0.2cm}\\
\displaystyle 0\neq X h^2(0,0) = Y^2f(0,0) =\frac{\partial^2 f}{\partial x^2}(0,0).
\end{array}
\end{equation}
In addition, $0$ is regular value of $f,$ thus we conclude that
  $$\frac{\partial f}{\partial y}(0) \neq 0. $$
Therefore, by the {\it Implicit Function Theorem}, there exists a unique $\mathcal{C}^{\omega}$  (resp. $\mathcal{C}^{\infty})$ function $\Phi: (-\varepsilon,\varepsilon)\to\R,$ defined for some small $\varepsilon>0,$ satisfying 
\begin{equation}\label{planarR2}
\Phi(0) =0\text{ and } f(x,\Phi(x))=0\ \text{for every} \ x\in(-\varepsilon,\varepsilon).
\end{equation}
This means that $U_1$ can be taken smaller, if necessary, in order that 
$$
\Sigma_1=\{(x,y)\in\R^2:\, y=\Phi(x), \, x\in (-\varepsilon,\varepsilon)\}.
$$

Now, we shall use Theorem \ref{thm:preparation} to construct a second transformation. From \eqref{planarR1} and \eqref{planarR2}, we compute
$$\Phi'(0)=0 \text{ and } \Phi''(0)=-\left(\dfrac{\partial f}{\partial y}(0,0)\right)^{-1}\dfrac{\partial ^2 f}{\partial x^2}(0,0)\neq0. $$
Defining $\tilde{\Phi}(x,y)\defeq \Phi(x),$ it is clear that
$$\tilde\Phi(0,0) =\frac{\partial\tilde\Phi}{\partial x}(0,0) = 0 \text{ and } \frac{\partial^2 \tilde\Phi}{\partial x^2}(0,0)\neq 0. $$
Thus, by Theorem \ref{thm:preparation}, there exists $\mathcal{C}^{\omega}$  (resp. $\mathcal{C}^{\infty})$ real-valued functions $\tilde a,\tilde b,$ defined in a neighborhood of $(0,0),$ such that 
$$x^2 + x\, \tilde a(\tilde\Phi(x,y),y) = \tilde b (\tilde \Phi(x,y),y).$$ 
Taking $a(x) \defeq  \tilde a(x,0)$ and $b(x)\defeq  \tilde b(x,0),$ we are able to conclude that
\begin{equation}\label{vida}
x^2 + x a(\Phi(x)) = b(\Phi(x)),
\end{equation}
for every $x$ in a neighborhood of $0.$  Clearly, $b(0)=0.$ Moreover, computing the first and second derivative of \eqref{vida} at $x=0$ and using \eqref{planarR2} we get, implicitly, that
$a(0)=0$ and $b'(0) \neq 0.$
Hence, consider the map 
\[
\beta(x,y)=\left(x +  \frac{a(y)}{2},  b(y) + \frac{a(y)^2}{4}\right),
\]
which is defined in a neighbourhood of $(0,0).$ Notice that $\det(D\beta(0,0))=b'(0)\neq0.$ Thus, there exists a neighborhood $U_2\subset V_1$ of $(0,0)$ on which $\beta$ is a diffeomorphism onto $W\defeq \beta(U_2).$ 
It is easy to see that $$Y^*(\tilde x,\tilde y)\defeq\beta_*Y(\tilde x,\tilde y)=(\beta\circ\alpha)_*X(\tilde x,\tilde y)=(1,0),\ \text{for every} \ (\tilde x,\tilde y)\in W.$$  Notice that, for $(\tilde x,\tilde y)=\beta(x,y),$ 
\begin{align*}
\tilde{y} - \tilde{x}^2 &= b(y) + \frac{a(y)^2}{4} - \left(x+\frac{a(y)}{2}\right)^2\\
&=  b(y) + \frac{a(y)^2}{4} - x^2 - x a(y) - \frac{a(y)^2}{4}\\
&= b(y) - a(y) x - x^2. 
\end{align*}
Thus, if $(x,y)\in\Sigma_1\cap U_2,$ that is, $y= \Phi(x),$ then the above identity together with \eqref{vida} imply that $\tilde y=\tilde x^2.$ This means that, in this new coordinate system, $\Sigma_1$ is transformed into $$\Sigma_2\defeq\beta(\Sigma_1\cap U_2)=\{(\tilde x,\tilde y)\in W:\, \tilde y=\tilde x^2\}.$$

As the last transformation, consider the following map,
$$\gamma(\tilde x,\tilde y)\defeq\left(\frac{1}{2}\left(\tilde x^2 - \tilde y\right),\tilde x\right), \ (\tilde x,\tilde y)\in W.$$
Notice that $\gamma$ is a local diffeomorphism around $(\tilde x,\tilde y)=(0,0).$ Therefore, we can take $W$ smaller, if necessary, in order that $\gamma$ is a diffeomorphism onto $V\defeq\gamma(W).$ It is easy to see that 
$$\gamma_*Y^*(x_1,x_2)=(\gamma\circ\beta\circ\alpha)_*X(x_1,x_2)=(x_2,1),\  \text{for every} \  (x_1,x_2) \in V $$
and
$$\gamma(\Sigma_2)=\{(x_1,x_2)\in V:\, x_1=0\}.$$
Finally, taking $U\defeq(\beta\circ\alpha)^{-1}(W)$ and  $\psi\defeq\gamma\circ\beta\circ\alpha\big|_U:U\rightarrow V,$ we get the following result proved.

\begin{theorem}
Let $X$ be a $\mathcal{C}^{\omega}$ (resp. $\mathcal{C}^{\infty}$) planar vector field and $\Sigma$ a  $\mathcal{C}^{\omega}$ (resp. $\mathcal{C}^{\infty}$) embedded planar curve.
Suppose that $p$ $\in\Sigma$ is a fold contact between $X$ and $\Sigma.$ Then, there exists a $\mathcal{C}^{\omega}$ (resp. $\mathcal{C}^{\infty}$) diffeomorphism  $\psi:U\rightarrow V,$ where $U,V\subset \R^2$ are respectively neighborhoods of $p$ and $0,$ such that $\psi(\Sigma\cap U) =\{(x_1,x_2)\in V:\, x_1=0\}$ and 
$\psi_*X(x_1,x_2) = (x_2, 1)$ for every $(x_1,x_2)\in V.$
\end{theorem}

\subsection{Proof of Theorem \ref{thm:vishik}}\label{proof:vishik}
Since $\Sigma$ is a codimension-1  $\mathcal{C}^{\omega}$ (resp. $\mathcal{C}^{\infty}$) embedded submanifold of $M,$ there exist a neighborhood $U_p\subset M$ of $p$ and a $\mathcal{C}^{\omega}$ (resp. $\mathcal{C}^{\infty}$) function $h: U_p\to \mathbb{R}$ such that $\Sigma\cap U_p=h^{-1}(0)$ and $\nabla h(\x)\neq0$ for every $\x\in\Sigma\cap U_p.$

We anticipate that, as in the planar case, the map $\psi:U\rightarrow V$ will be given as the composition of three maps,  $\psi=\gamma\circ\beta\circ\alpha.$ In what follows, we shall construct each one of these maps.

The first map $\alpha$ will be given as the composition of two maps, $\alpha=\check\alpha\circ \hat \alpha,$ and has the purpose of transforming the boundary equation in such way that we can use Theorem \ref{thm:preparation}. As in the planar case, the first one $\hat\alpha$ is obtained through the {\it Tubular Flow Theorem}, the second one $\check \alpha$ is an auxiliary map.  Since $X(p)\neq0,$ the {\it Tubular Flow Theorem} provides a  $\mathcal{C}^{\omega}$ (resp. $\mathcal{C}^{\infty})$ chart $(\hat \alpha,U_1)$ of $p$ satisfying  
$$\widehat Y(\hat \y)\defeq\hat \alpha_*X(\hat \y) = (1,0,\ldots,0),\ \text{for every} \ \hat \y \in \widehat V_1\defeq \hat \alpha(U_1)\subset\R^m.$$ 
Defining $\hat f = h\circ \hat \alpha^{-1},$ we have that $\Sigma$ is transformed into $\widehat \Sigma_1\defeq\hat \alpha(\Sigma\cap U_1)=\hat f^{-1}(0).$ In addition, from \ref{prop21}, we have that 
$$0=Xh(p) = Y \hat f(0)=\frac{\partial \hat f}{\partial \hat y_1}(0).$$
Here, we are denoting $\hat \y=(\hat y_1,\ldots,\hat y_m).$ Since $\nabla \hat f(\hat \y)\neq0$ for every $\hat \y \in\widehat\Sigma_1,$ we can assume, without loss of generality, that
$$\frac{\partial \hat f}{\partial \hat y_m}(0) \neq 0.$$
Therefore, by the {\it Implicit Function Theorem}, there exists a unique $\mathcal{C}^{\omega}$  (resp. $\mathcal{C}^{\infty})$ function $\Phi: B_{\varepsilon}^{m-1}\subset\R^{m-1}\to\R,$ defined for some  $\varepsilon>0,$ satisfying $\Phi(0) =0$ and
$$
\hat f(\hat y_1,\ldots \hat y_{m-1},\Phi(\hat y_1,\ldots \hat y_{m-1}))=0,\ \text{for every} \ (\hat y_1,\ldots \hat y_{m-1})\in B_{\varepsilon}^{m-1}.
$$
Here,  $B_{\varepsilon}^{m-1}$ denotes the open ball of radius $\varepsilon$ centered at the origin of $\R^{m-1}.$ This means that $U_1$ can be taken smaller, if necessary, in order that 
$$
\widehat \Sigma_1=\{\hat \y\in\R^m:\, \hat y_m=\Phi(\hat y_1,\ldots,\hat  y_{m-1}), \, (\hat y_1,\ldots, \hat y_{m-1})\in  B_{\varepsilon}^{m-1}\}.
$$
In the planar case, Theorem \ref{thm:preparation} was applied directly for the function $\Phi$. Here, before using Theorem \ref{thm:preparation}, we have to consider the following auxiliary diffeomorphism 
$$
\check\alpha (\hat y_1,\ldots,\hat y_m)\defeq\left(\hat y_1\ , \ \ldots \ ,\ \hat y_{m-1},\  \hat y_m - \sum_{i=2}^{m-1} \frac{\partial \Phi}{\partial \hat y_i}(0) \hat y_i\right),\ (\hat y_1,\ldots,\hat y_m)\in \widehat V_1.
$$

Denote $\y=( y_1,\ldots, y_m)$. It is easy to see that
$$ Y( \y)\defeq\check\alpha_*\widehat Y( \y) = (1,0,\ldots,0),\ \text{for every} \  \y \in V_1\defeq \check\alpha(\widehat V_1)\subset\R^m,$$ 
and
 $$ \Sigma_1\defeq \check\alpha(\widehat \Sigma_1)=\{ \y\in  V_1:\,  y_m=\phi( y_1,\ldots, y_{m-1}),\, ( y_1,\ldots  y_{m-1})\in B_{\varepsilon}^{m-1}\},$$
where
$$
\phi ( y_1,\ldots  y_{m-1}) \defeq \Phi( y_1,\ldots, y_{m-1})-\sum_{i=2}^{m-1} \frac{\partial \Phi}{\partial \hat y_i}(0)  y_i.
$$
Notice that 
\begin{equation}\label{dphi}
\frac{\partial\phi}{\partial y_i}(0) = 0, \text{ for } \ i\in \{1,\ldots,m-1\}.
\end{equation}
In addition, denoting $\alpha\defeq\check\alpha\circ\hat\alpha$ and $f\defeq h\circ \alpha^{-1},$ we get
\begin{equation}\label{generalR2}
 f( y_1,\ldots  y_{m-1},\phi( y_1,\ldots  y_{m-1}))=0,\ \text{for every} \ ( y_1,\ldots  y_{m-1})\in B_{\varepsilon}^{m-1}.
\end{equation}

Now, in order to construct the second transformation $\beta$, we shall apply Theorem \ref{thm:preparation} for the function $\phi$. From  \ref{prop21}, we have that
\begin{equation}\label{generalR0}
X^i h(\y) =  Y^i  f(\alpha (\y)) =\frac{\partial^i  f}{\partial  y_1^i}(\alpha(\y)).
\end{equation}
Therefore,
\begin{equation}\label{generalR1}
\begin{array}{l}
\displaystyle 0=X^i h(p) =  Y^i  f(0) =\frac{\partial^i  f}{\partial  y_1^i}(0),\ \text{ for } \ i\in\{1,2,\ldots,k\},\text{ and }\vspace{0.2cm}\\
\displaystyle 0\neq X^{k+1} h(p) =  Y^{k+1}  f(0) =\frac{\partial^{k+1}  f}{\partial  y_1^{k+1}}(0).
\end{array}
\end{equation}
Thus, from \eqref{generalR1} and \eqref{generalR2}, we compute  
\begin{equation}\label{dkphi}
\begin{array}{l}
\displaystyle\frac{\partial^i \phi}{\partial y_1^i}(0) = -\left(\dfrac{\partial f}{\partial y_m}(0)\right)^{-1}\dfrac{\partial ^i f}{\partial y_1^i}(0)=0, \text{ for }1\leq i \leq k,\text{ and }\vspace{0.2cm}\\
\displaystyle \frac{\partial^{k+1} \phi}{\partial y_1^{k+1}}(0) = -\left(\dfrac{\partial f}{\partial y_m}(0)\right)^{-1}\dfrac{\partial ^{k+1} f}{\partial y_1^{k+1}}(0)\neq0.
 \end{array}
\end{equation}
Now, by Theorem \ref{thm:preparation}, there exists $\mathcal{C}^{\omega}$  (resp. $\mathcal{C}^{\infty})$ real-valued functions $a_1,\ldots,a_k,b:\mathbb{R}^{m-1} \to \mathbb{R},$ defined for $(y_2,\ldots,y_m)$ in a neighborhood of $0\in\R^{m-1},$ such that
\begin{equation}
    y_1^{k+1} + \sum_{i=1}^{k} y_1^i\, a_i(y_2,\ldots,y_{m-1},\phi(y_1,\ldots,y_{m-1})) =  b(y_2,\ldots,y_{m-1},\phi(y_1,\ldots,y_{m-1} ))\label{mc2},
\end{equation} 
Clearly, $b(0)=0.$ Computing the derivative of \eqref{mc2} in the variable $y_i,$ for $2\leq i \leq m-1,$ at $\y=0$ and using \eqref{dphi}  we get
\begin{equation}\label{dbdyi}
\frac{\partial  b}{\partial y_i}(0) = 0,\text{ for }2\leq i \leq m-1. 
\end{equation}
Moreover, computing the $i$th-derivative of \eqref{mc2} in the variable $y_1$ at $\y=0$ for $1\leq i\leq k$ and $i=k+1$ we obtain
\begin{equation}\label{ai}
a_i(0)=0,\text{ for } 1\leq i\leq k,
\end{equation}
and
\begin{equation}\label{partialky1}
 \frac{\partial }{\partial y_1^{k+1}}\Big(b(y_2,\ldots,y_{m-1},\varphi(y_1,\ldots,y_{m-1}))\Big)\Big|_{\y=0}=(k+1)!\,,
 \end{equation}
 respectively. Thus, from \eqref{dphi} and \eqref{partialky1}, we conclude that
\begin{equation}\label{db} \frac{\partial b}{\partial y_m}(0) \neq 0.
\end{equation}
Now, we claim that
 the matrix
\begin{align*}
A\defeq \left(\frac{\partial a_i}{\partial y_j}(0)\right)_{(i,j)\in \{1,\ldots,k-1\}\times\{2,\ldots,k\}}
\end{align*}
is invertible. Notice that, applying $\displaystyle\frac{\partial^{i+1}}{\partial y_1^i \partial y_j}\Big|_{\y=0}$ to both sides of the equation \eqref{mc2}, we get that
$$\frac{\partial a_i}{\partial y_j}(0) = \frac{\partial b}{\partial  y_m}(0)\cdot \frac{\partial^{1+i}\phi }{\partial y_1^i\partial y_j}(0).$$
Thus, from \eqref{db}, $A$ is invertible if, and only if, 
$$B\defeq\left(\dfrac{\partial^{1+i}\phi}{\partial y_1^i\partial y_j}(0)\right)_{(i,j)\in \{2,\ldots,k\}\times\{2,\ldots,k\}}$$
is invertible. Denote
$$\overline{B}\defeq \left(\frac{\partial^{i+1} \phi}{\partial y_1^i \partial y_j}(0) \right)_{(i,j)\in \{1,\ldots,k\}\times\{1,\ldots,k\}}.$$
From \eqref{dkphi},
$$
\overline{B} = \left(\begin{array}{c|ccc}
    0 &         &   &       \\ 
    \vdots &  \multicolumn{3}{c}{B} \\ 
    0 &  &     & \\ 
    \hline
    \displaystyle\frac{\partial^{k+1}\phi}{\partial y_1^{k+1}}(0) &    * &\ldots &   * 
    \end{array}\right).
$$
Since $ \displaystyle\frac{\partial^{k+1}\phi}{\partial y_1^{k+1}}(0)\neq 0,$  we conclude that $B$ is invertible if, and only if, $\overline{B}$ is invertible. Thus, in what follows, we shall prove that $\overline{B}$ is invertible. 
By hypothesis, we know that $$\{\nabla h(p), \nabla Xh(p), \ldots, \nabla X^k h(p)\}$$is linearly independent. Thus, from \eqref{generalR0}, we get that 
\begin{equation}\label{li}
\left\{\nabla f, \nabla \frac{\partial f}{\partial y_1}(0), \nabla \frac{\partial^2 f}{\partial y_1^2}(0), \ldots, \nabla \frac{\partial^k f}{\partial y_1^{k}}(0)  \right\}
\end{equation} is also linearly independent. Since the definition of simple $k$-contact does not depend on the function $f$, which describes the manifold $\Sigma_1$ as an inverse image of regular value, by \eqref{generalR2} we can assume in \eqref{li} that $f=y_m-\phi(y_1,\ldots,y_{m-1}).$ Thus,  we conclude that the matrix
$$ \left(\frac{\partial^{i+1} \phi}{\partial y_1^i \partial y_j} \right)_{(i,j)\in \{1,\ldots,k\}\times\{1,\ldots,m-1\}}$$
has rank $k.$  
Therefore, up to a permutation of the coordinates, we can assume, without loss of generality, that the matrix $\overline{B}$ is invertible and, therefore, we conclude that the matrix $A$ is invertible. Hence, consider the map $\beta(\y)=\big(\beta_1(\y),\ldots,\beta_m(\y)\big)$ defined for $\y\in V_1$ by
$$\begin{array}{l}
        \beta_1(\y) = y_1 +\dfrac{ a_k(y_2,\ldots,y_m)}{k+1}, \vspace{0.3cm}\\
        \!\!\!\begin{array}{ll}
        \beta_{k+1-j}(\y) =& \displaystyle\binom{k+1}{j}  \left(\frac{-a_k(y_2,\ldots,y_m)}{k+1}\right)^{k+1-j} \vspace{0.2cm} \\
      & + \displaystyle \sum_{i=j}^{k} \binom{i}{j} \left(\frac{-a_k(y_2,\ldots,y_m)}{k+1}\right)^{i-j} a_i(y_2,\ldots,y_m), \text{ for } 1\leq j \leq k-1,
       \end{array}\vspace{0.3cm}\\
    
        \!\!\!\begin{array}{ll}\beta_{k+1}(\y) =&\displaystyle -b(y_2,\ldots,y_m) + k\left(\frac{-a_k(y_2,\ldots,y_m)}{k+1}\right)^{k+1} \vspace{0.2cm}\\ &\displaystyle+ \sum_{i=1}^{k-1}\left(\frac{-a_k(y_2,\ldots,y_m)}{k+1}\right)^ia_i(y_2,\ldots, y_m),\end{array}\vspace{0.3cm}\\
        
        \beta_{i}(\y)=y_i, \,\, \text{ for }\ k+2 \leq i \leq m-1,  \vspace{0.3cm} \\
        
       \beta_m(\y) = y_{k+1}.    
        \end{array} $$

From \eqref{dbdyi},  \eqref{ai}, and \eqref{db}, we get that $$\det(D\beta(0))=-\dfrac{\partial b}{\partial y_m}(0)\det(A)\neq0.$$

Thus, there exists a neighborhood $U_2\subset V_1$ of $0$ on which $\beta$ is a diffeomorphism onto $W\defeq \beta(U_2).$ The motivation for this choice of coordinates will become clear further ahead. Notice that
$$Y^*(\tilde \x)\defeq\beta_*Y(\tilde \x)=(\beta\circ\alpha)_*X(\tilde \x)=(1,0,\ldots,0),\ \text{for every} \ \tilde \x\in W.$$ 
Denote $\tilde\x=(\tilde x_1,\ldots,\tilde x_m).$ We claim that 
\begin{equation}\label{CLAIM1}
\Sigma_2\defeq \beta(\Sigma_1\cap U_2)=\{ \tilde \x\in  W:\,  \tilde x_1^{k+1} + \tilde x_2 \tilde x_1^{k-1} + \tilde x_3 \tilde x_1^{k-2} + \ldots + \tilde x_k \tilde x_1 + \tilde x_{k+1}=0\}.
\end{equation}
Indeed, take $\y=(y_1,\ldots,y_m)\in\Sigma_1.$ We know that
$$y_{m}= \phi(y_1,\ldots,y_{m-1}).$$
Thus, from \eqref{mc2}, the following equation is fulfilled
\begin{equation}
y^{k+1} + \sum_{i=1}^{k} y_1^i\cdot a_i(y_2,\ldots,y_{m}) =  b(y_2,\ldots,y_{m}).\label{cscovarde}
\end{equation}
For the sake of simplicity, denote $a_i\defeq a_i(y_2,\ldots,y_m)$ and $b\defeq b(y_2,\ldots,y_m).$ Thus, 
\begin{align*}
0&=y_1^{k+1} + \sum_{i=1}^{k} y_1^i\cdot a_i -  b\\
&= \left(y_1  + \frac{a_k}{k+1}\right)^{k+1} - \sum_{i=0}^{k}\binom{k+1}{i}y^i_1\left(\frac{a_k}{k+1}\right)^{k+1-i} + a_k y_1 + \sum_{i=1}^{k-1} y_1^{i}a_i - b\\
&=\left(y_1  + \frac{a_k}{k+1}\right)^{k+1} - \sum_{i=0}^{k-1}\binom{k+1}{i}y^i_1\left(\frac{a_k}{k+1}\right)^{k+1-i}  + \sum_{i=1}^{k-1} y_1^{i}a_i - b.
\end{align*}
Now, take $\tilde \x=(\tilde x_1,\ldots,\tilde x_m)=\beta(\y).$ In particular, $\tilde{x}_1=y_1 + a_k/(k+1).$ Hence,
\begin{align}
0=&\tilde{x}_1^{k+1} - \sum_{i=0}^{k-1}\binom{k+1}{i} \left(\tilde{x}_1 - \frac{a_k}{k+1}\right)^i\left(\frac{a_k}{k+1}\right)^{k+1-i}  + \sum_{i=1}^{k-1} \left(\tilde{x}_1  - \frac{a_k}{k+1}\right)^{i}a_i - b\nonumber\\
=&\tilde{x}_1^{k+1} - \sum_{i=0}^{k-1} \sum_{j=0}^i\binom{k+1}{i} \binom{i}{j}(-1)^{i-j}\tilde{x}_1^j  \left(\frac{a_k}{k+1}\right)^{k+1-j}  + \sum_{i=1}^{k-1}\sum_{j=0}^i \binom{i}{j}\tilde{x}_1^j \left(\frac{-a_k}{k+1}\right)^{i-j} a_i - b\nonumber\\
=& \tilde{x}_1^{k+1} - \sum_{j=1}^{k-1} \sum_{i=j}^{k-1}\binom{k+1}{i} \binom{i}{j}(-1)^{i-j}\tilde{x}_1^j  \left(\frac{a_k}{k+1}\right)^{k+1-j}  + \sum_{j=1}^{k-1}\sum_{i=j}^{k-1} \binom{i}{j}\tilde{x}_1^j \left(\frac{-a_k}{k+1}\right)^{i-j} a_i \nonumber\\
&-b +  k  \left(\frac{-a_k}{k+1}\right)^{k+1} + \sum_{i=1}^{k-1}  \left(\frac{-a_k}{k+1}\right)^{i} a_i\nonumber \\
=& \tilde{x}_1^{k+1} + \sum_{j=1}^{k-1} \tilde{x}_1^j \left(\sum_{i=j}^{k-1}\binom{k+1}{i} \binom{i}{j}(-1)^{i+1-j} \left(\frac{a_k}{k+1}\right)^{k+1-j}  + \sum_{i=j}^{k-1} \binom{i}{j} \left(\frac{-a_k}{k+1}\right)^{i-j} a_i\right)\nonumber\\
&-b +  k  \left(\frac{-a_k}{k+1}\right)^{k+1} + \sum_{i=1}^{k-1}  \left(\frac{-a_k}{k+1}\right)^{i} a_i\nonumber  \\
=& \tilde{x}_1^{k+1} + \sum_{j=1}^{k-1} \tilde{x}_1^j \left(\sum_{i=j}^{k}\binom{k+1}{i} \binom{i}{j}(-1)^{i+1-j} \left(\frac{a_k}{k+1}\right)^{k+1-j}+  \sum_{i=j}^{k} \binom{i}{j} \left(\frac{-a_k}{k+1}\right)^{i-j} a_i\right)\nonumber\\
&-b +  k  \left(\frac{-a_k}{k+1}\right)^{k+1}+ \sum_{i=1}^{k-1}  \left(\frac{-a_k}{k+1}\right)^{i} a_i \nonumber  \\
=& \tilde{x}_1^{k+1} + \sum_{j=1}^{k-1} \tilde{x}_1^j \left(\binom{k+1}{j}  \left(\frac{-a_k}{k+1}\right)^{k+1-j} +  \sum_{i=j}^{k} \binom{i}{j} \left(\frac{-a_k}{k+1}\right)^{i-j} a_i\right)- b + k  \left(\frac{-a_k}{k+1}\right)^{k+1} \label{milagre} \\
&+\sum_{i=1}^{k-1}  \left(\frac{-a_k}{k+1}\right)^{i} a_i \nonumber \\
=&\tilde x_1^{k+1} + \tilde x_2 \tilde x_1^{k-1} + \tilde x_3 \tilde x_1^{k-2} + \ldots + \tilde x_k \tilde x_1 + \tilde x_{k+1}.\nonumber
\end{align}
 Therefore, $\beta(\Sigma_1 \cap U_2) \subset \{ \tilde \x\in  W:\,  \tilde x_1^{k+1} + \tilde x_2 \tilde x_1^{k-1} + \tilde x_3 \tilde x_1^{k-2} + \ldots + \tilde x_k \tilde x_1 + \tilde x_{k+1}=0\}$. The opposite continence is obtained by noticing that the two previous sets are    codimension-1 embedded submanifolds of $W$. Thus, since $0$ lies in both submanifolds, by shrinking $W$, if necessary, we get that the sets coincide.
%
This proves (\ref{CLAIM1}). 

Notice that the above manipulation motivates the definition of $\beta$. Indeed, first we have defined $\beta_1(\y) = y_1 + a_{k}(y_2,\ldots,y_m)/(k+1)$ to get rid off the term $y_1^k$ on equation (\ref{cscovarde}). Then, the remaining $\beta_i$, $1<i\leq m$, are chosen from manipulation acquired in  (\ref{milagre}).

For the last transformation, consider the map $\gamma(\tilde\x)=\big(\gamma_1(\tilde\x),\ldots,\gamma_m(\tilde\x)\big)$ defined for $\tilde\x\in W$ by
$$   \begin{array}{l}
        \gamma_i(\tilde\x)=  \displaystyle\dfrac{\tilde x_1^{k+2-i}}{(k+2-i)!} + \sum_{l=2}^{k+2-i}\frac{(k+1-l)!}{(k+1)!} \frac{\tilde x_1^{k+2-i-l}\tilde x_i }{(k+2-i-l)!}, \text{ for } \ 1\leq i \leq k+1,\vspace{0.2cm}  \\
        \gamma_i(\tilde\x) = \tilde x_i, \text{ for } \ k+2\leq i \leq m.   \\   
        \end{array} $$
Notice that $\gamma$ is a local diffeomorphism around $\tilde\x=0.$ Therefore, we can take $W$ smaller, if necessary, in order that $\gamma$ is a diffeomorphism onto $V\defeq\gamma(W).$ Take, $\x=(x_1,\ldots x_m)$.  Then, for every $\x \in V,$
\begin{align*}
\gamma_*Y^*(\gamma(\x))&=\gamma_*(1,0,\ldots ,0)(\x)\\
&=\mathrm{d}\gamma_{\gamma^{-1}(\x)}(1,0,\ldots ,0)\\
&=\left(\frac{\partial\gamma_1}{\partial x_1}\left(\gamma^{-1}(\x)\right),\ldots,\frac{\partial\gamma_m}{\partial x_1}\left(\gamma^{-1}(\x)\right) \right).
\end{align*}
Let $\pi_i: \mathbb{R}^n \to \mathbb{R}^m$  be the projection of the $i$-th coordinate. Thus, for $1\leq i <k+1,$ we have
\begin{align*}
\pi_i\circ \gamma_*Y^*(\gamma(\x)) &= \frac{\partial \gamma_i}{\partial x_1}(\gamma^{-1}(\x))\\
&=  \left.\frac{\partial}{\partial \tilde{x}_1}\left(\displaystyle\dfrac{\tilde x_1^{k+2-i}}{(k+2-i)!} + \sum_{l=2}^{k+2-i}\frac{(k+1-l)!}{(k+1)!} \frac{\tilde x_1^{k+2-i-l}\tilde x_i }{(k+2-i-l)!}\right)\right|_{\tilde\x = \gamma^{-1}(\x)}\\
&=\left.\left(\displaystyle\dfrac{\tilde x_1^{k+2-(i+1)}}{(k+2-(i+1))!} + \sum_{l=2}^{k+2-(i+1)}\frac{(k+1-l)!}{(k+1)!} \frac{\tilde x_1^{k+2-(i+1)-l}\tilde x_i }{(k+2-(i+1)-l)!}\right)\right|_{\tilde\x = \gamma^{-1}(\x)}\\
&= \gamma_{i+1}(\tilde \x)\Big\vert_{\tilde\x = \gamma^{-1}(\x)} 
=\gamma_{i+1}\left(\gamma^{-1}(\x)\right) = x_{i+1}.
\end{align*}
For $i=k+1$, we have
\[
\pi_{k+1}\circ \gamma_*Y^*(\gamma(\x)) = \frac{\partial \gamma_k}{\partial x_1}(\gamma^{-1}(\x))
=  \left.\frac{\partial \tilde{x}_1}{\partial \tilde{x}_1}\right|_{\tilde\x = \gamma^{-1}(\x)}
=1.
\]
Finally, when $k+1< i\leq m,$ we have
\[
\pi_i\circ \gamma_*Y^*(\gamma(\x)) = \frac{\partial \gamma_i}{\partial x_1}(\gamma^{-1}(\x))
=  \left.\frac{\partial \tilde{x}_i}{\partial \tilde{x}_1}\right|_{\tilde\x = \gamma^{-1}(\x)}
=0.
\]
Therefore, we achieved
$$\gamma_*Y^*(\x)=(\gamma\circ\beta\circ\alpha)_*X(\x)=(x_2,x_3,\ldots,x_k,1,0\ldots,0), \text{ if } m>k+1$$
and
$$\gamma_*Y^*(\x)=(\gamma\circ\beta\circ\alpha)_*X(\x)=(x_2,x_3,\ldots,x_k,1), \text{ if } m=k+1.$$
In both cases,
\begin{align*}
\gamma(\Sigma_2)&=\gamma\left(\left\{ \tilde \x\in  W:\,  \tilde x_1^{k+1} + \tilde x_2 \tilde x_1^{k-1} + \tilde x_3 \tilde x_1^{k-2} + \ldots + \tilde x_k \tilde x_1 + \tilde x_{k+1}=0\right\}\right)\\
&= \gamma\left(\left\{ \tilde \x\in  W:\,  (k+1)!\cdot \gamma_1(\tilde\x)=0\right\}\right)\\
&= \gamma \circ \gamma_1^{-1}(0)\\
&=\{\x \in V;\ x_1=0\}.
\end{align*}

Finally, we conclude the proof of Theorem \ref{thm:vishik} by taking $U\defeq(\beta\circ\alpha)^{-1}(W)$ and  $\psi\defeq\gamma\circ\beta\circ\alpha\big|_U:U\rightarrow V.$

\section*{Acknowledgments}

The authors are very grateful to Marco A. Teixeira for meaningful discussions and constructive criticism on the manuscript.

MM is partially supported by FAPESP grants 2017/23692-6 and 2019/06873-2. RMM is partially supported by a FAPESP grant 2018/03338-6. DDN is partially supported by FAPESP grants 2018/16430-8 and 2019/10269-3, and by CNPq grants 306649/2018-7 and 438975/2018-9. RMM and DDN are also partially supported by a FAPESP grant 2018/13481-0.

 \bibliographystyle{plain}
\bibliography{bibliografia}
\end{document}